\documentclass{amsart}

\usepackage{hyperref}
\usepackage{amsmath,color,graphicx,amssymb,amscd}
\usepackage{degt} 
\usepackage{tikz}
\usepackage{tikz-cd}

\usepackage{booktabs}
\usepackage{comment} 
\usepackage[paper=portrait,pagesize]{typearea}
\usepackage{longtable}
\usepackage{array, boldline, makecell}

\newlength\mylength

\def\L{\mathbf{L}}

\def\M{\mathcal{M}}

\def\rk{\mathrm{rk}}
\def\discr{\operatorname{discr}}

\def\Aut{\operatorname{Aut}}

\def\ie{\emph{i.e\PERIOD}}

\def\iq.{\emph{i.q\PERIOD}}
\def\cf.{\emph{cf.}}
\def\viz.{\emph{viz\PERIOD}}
\def\vs.{\emph{vs\PERIOD}}

\begin{document}

\title[Projective models of $K3$--surfaces]{Classification of Projective models of $K3$--surfaces}

\author[\c{C}.~G\"{u}ne\c{s} ~Akta\c{s} ]{\c{C}\.{i}sem G\"{u}ne\c{s} Akta\c{s}}
\address{
 Department of Engineering Sciences, Abdullah G\"{u}l University\\
 38080, Kayseri, Turkey} \email{cisem.gunesaktas@agu.edu.tr}

\thanks{The author was supported by the \textit{T\"{U}B\.{I}TAK} grant $123$F$111$}

\subjclass[2010]{Primary 14J28; Secondary  14J10, 14J17}

\keywords{$K3$-surface, projective model, singular model, simple singularity}

\date{}

\dedicatory{}

\begin{abstract}

Extending the methods previously established for plane sextics  and 
quartic surfaces, we develop a unified and systematic framework for the equisingular deformation classification of projective models of $K3$-surfaces of arbitrary polarization. We treat all polarizations of degree $2n$ with 
$n\le 4$, including the higher-degree cases of sextic surfaces in $\mathbb{P}^{4}$ 
and octic surfaces in $\mathbb{P}^{5}$.
 
The computational machinery developed here is degree-independent and provides 
a complete tool for enumerating and distinguishing the deformation classes 
arising from singular projective models of $K3$-surfaces. This framework lays the foundation for future investigations in higher degrees, 
particularly the geometry of equisingular strata, extremal configurations, 
real structures.

\end{abstract}


\maketitle


\section{Introduction}

K3-surfaces play a central role in algebraic geometry, bridging the gap between geometric complexity and their rigorous lattice structure. While their underlying
geometry is subtle and poses challenging problems, they also offer
powerful tools, such as the global Torelli theorem~\cite{K3}  and the surjectivity of the period map~\cite{periodmap} to address these problems. These foundational results
allow one to translate geometric questions about K3-surfaces into purely
arithmetical ones involving their N\'{e}ron--Severi lattices, which can then
be studied by using Nikulin's~\cite{Niku2} theory of discriminant forms combined with the theory of Miranda and Morrison~\cite{MM1,MM2,MM3}.

A projective model $X \to \mathbb{P}^{n+1}$ of a K3-surface $X$ gives rise
to a polarization, \ie, a class $h \in H_2(X)$ of the hyperplane section,
so that $h^2 = 2n$. The cases $h^2=2$ (sextic curves in $\mathbb{P}^2$)
and $h^2=4$ (quartic surfaces in $\mathbb{P}^3$) have been extensively
studied. By translating geometric questions into arithmetic problems via
the N\'{e}ron--Severi lattice and discriminant form theory, complete equisingular classifications have emerged. Using the
arithmetical reduction~\cite{Alex1}, Shimada~\cite{Shimada.Maximizing} gave a complete
description of the moduli spaces of maximizing sextics. Later, based on
the same approach, Degtyarev and Akyol~\cite{Alex2} completed the equisingular
deformation classification of plane sextics. This approach was also used
by G\"{u}ne\c{s} Akta\c{s}~\cite{Cisem1,aktacs2019real,Cisem2024} to obtain a complete description of the
equisingular strata of the simple quartics. In the meanwhile, Shimada~\cite{Shimada.connEllK3} has listed the connected components of the moduli space of the
Jacobian elliptic K3-surfaces.

The study of singular projective models of K3-surfaces is a vibrant area of
research, with significant progress made for low polarizations
($h^2=2,4$) and promising avenues for higher polarizations ($h^2=6,8$)
which remain underexplored despite their theoretical promise. The
well-established results and techniques obtained for the former cases can
serve as a foundational framework for advancing the equisingular
deformation classification of the latter. The principal novelty of this paper
is to extend the methods originally originally developed
by Degtyarev–Akyol~\cite{Alex2} for sextic curves, and by G\"{u}ne\c{s} Akta\c{s}~\cite{Cisem1,aktacs2019real,Cisem2024} for quartic surfaces, to all projective models of K3-surfaces of arbitrary degrees yielding a unified computer-aided (\texttt{GAP}~\cite{GAP}) framework for the equisingular deformation classification of $K3$-surfaces with any polarization. This study hinges on three pillars:
\begin{itemize}
	\item \textbf{Arithmetic Reduction:} Translating geometric complexity into
	lattice-theoretic problems via N\'{e}ron--Severi lattices.
	\item \textbf{Computational Rigor:} Employing \texttt{GAP} (2015) to
	manage combinatorial explosions in higher polarizations, adapting
	existing codes for sextic curves/quartic surfaces to sextic surfaces and
	octic surfaces.
	\item \textbf{Synthesis of Techniques:} Integrating Nikulin's discriminant
	forms and Miranda--Morrison’s theory.
\end{itemize}
From a theoretical standpoint, higher polarizations ($h^2=6$, $h^2=8$) offer
novel challenges due to increased dimensionality and interaction between
polarization degree and Picard rank. Yet, they also provide fertile ground by
representing natural extensions. These models, while less studied, promise
insights into the interplay between polarization degree, singularity types,
and lattice structures. By extending the theory to higher polarizations,
this work aims to unify fragmented classifications into a universal
framework, ultimately contributing to a deeper understanding of K3-surface
moduli and their geometric realizations. This project not only bridges gaps
in equisingular classification but also opens avenues for exploring real
algebraic geometry and the computational limits of lattice-based
algorithms in higher dimensions.

\subsection{History of the problem}

One of the major questions in singularity theory is the extent to which the
local data about a singular object (say, types of singularities) determine
the object globally, up to appropriately defined deformation. Various
deformation classification problems for K3-surfaces have been intensively
studied in the literature. Any model of K3-surfaces is necessarily simple,
i.e., has $\mathbf{A}$-$\mathbf{D}$-$\mathbf{E}$ type simple singularities. Irreducible smooth
rational curves on $X$ are (some) classes $\ell$ in the N\'{e}ron--Severi
lattice $\mathrm{NS}(X)\subset H_2(X)$ of square $-2$ intersecting $h$ in a
prescribed way. The simplest case is that of exceptional divisors,
$\ell\cdot h=0$, in which we are speaking about the equisingular deformation
classification of singular models. Singular sextic curves and quartic
surfaces are a special case of polarized K3-surfaces with prescribed
combinations of exceptional divisors. More generally, given a K3-surface
$X$ and a polarization $h\in\mathrm{NS}(X)$, a smooth rational curve
$L\subset X$ realizing a class $\ell := [L]\in\mathrm{NS}(X)$ is
\begin{itemize}
	\item an exceptional divisor, if $\ell\cdot h = 0$,
	\item a line, if $\ell\cdot h = 1$,
	\item a conic, if $\ell\cdot h = 2$,
	\item a twisted cubic, if $\ell\cdot h = 3$,
\end{itemize}
etc. By the Hodge index theorem, the number of classes
$\ell \in \mathrm{NS}(X)$, $\ell^2 = -2$, $\ell \cdot h \leq \text{const}$
is always finite, and one can naturally consider configurations of rational
curves of any of these kinds or even combinations thereof; however,
currently, only configurations of lines seem to be within reach.


Singular models of K3-surfaces can be regarded as a variation of the line
counting problem: instead of lines, we consider rational curves orthogonal
to the polarization (exceptional divisors). As
a very special case of polarizations, namely, $h^2=0$, one can mention
Shimada~\cite{Shimada.connEllK3} dealing with Jacobian elliptic
K3-surfaces. The deformation classification
of plane algebraic curves, either real or singular, is also a well-established
subject. The classification of cubic and quartic curves goes as far back as
to, respectively, Newton~\cite{Newton1967} and Klein~\cite{Klein1922}. After the classification of singular quintics was completed in Degtyarev~\cite{Degtyarev1989a}, sextics have
naturally become an object of intensive study. The na\"{i}ve equation-based
approach works well for sextics with a non-simple singular point (see
Degtyarev~\cite{Degtyarev1989NonSimple}. The same approach for simple sextics is exploited in
Artal et~al.~\cite{Artal2001,Artal2002,Artal2007}, Oka and Po~\cite{OkaPo2002} and a great number of other papers
(see Degtyarev~\cite{Degtyarev2011} for a recent survey and more details). Unfortunately,
this direction is bound to fail: even when the classification is known, it requires a tremendous amount of work to find an explicit equation of any given curve.

The K3-theoretical approach was used by Persson~\cite{Persson}, who established the bound on the total Milnor number, and Urabe~\cite{Urabe1988,Urabe1989} and Yang~\cite{Yang.sextics},
who gave a classification of the possible sets of singularities. However,
it was not until Degtyarev~\cite{Alex1}, where the arithmetical reduction
theorem was stated and a few examples were considered, that it was realized
that this approach is the key to the complete classification. This was
followed by a number of papers, see, e.g., Degtyarev~\cite{Degtyarev2008b},
Shimada~\cite{Shimada:Zsplitting}; the classification was finally settled for the
so-called maximizing curves in Shimada\cite{Shimada.Maximizing}; and in general in Akyol and Degtyarev~\cite{Alex1}.

The next natural application of the theory is to simple quartic surfaces,
$h^2=4$. Quartics play the same role in the realm of spatial surfaces as
sextics do for curves, thus can be treated by
similar methods. Surfaces of degree $3$ were classified, by the classical
tools, in Schläfli~\cite{Schlafli1863} and Klein~\cite{Klein1873}. Surfaces of degree $4$ with a
non-simple singular point can still be handled classically (see Degtyarev~\cite{Degtyarev1989,Degtyarev1994II}, but no attempts seem to have been made to extend this approach
to surfaces with simple singularities only, apparently because of the great
complexity of the problem. The K3-theoretical approach was applied
successfully to the classification of smooth real quartics in Kharlamov~\cite{Kharlamov}
 and Nikulin~\cite{Niku1}; a few sporadic examples of the study
of a certain special class of singular surfaces related to semidefinite
programming can be found in Degtyarev and Itenberg~\cite{DegtyarevItenberg2011} and Ottem
et~al.~\cite{Ottemetal}. The classification of the sets of simple singularities
realized by quartic surfaces is found in Urabe~\cite{Urabe3,Urabe2,Urabe1} and Yang~\cite{Yang}. The systematic classification of simple singular quartic surfaces
was started by G\"{u}ne\c{s} Akta\c{s}~\cite{Cisem1}, where she completed the classification
of non-special quartics, then it was completed in G\"{u}ne\c{s} Akta\c{s}~\cite{aktacs2019real,Cisem2024},
where the ultimate answer was given in the spirit of Akyol and Degtyarev~\cite{Alex2}.


Hypersurfaces of degree $6$ in $\mathbb{P}^4$ represent the first
non-trivial case of K3-surfaces with polarization $h^2=6$. Their geometry
is inherently more complex due to the higher-dimensional ambient space and
the interplay between $h$ and the Picard lattice. While Beauville~\cite{Beauville1983}
established foundational results on K3-hypersurfaces, the classification of
their singular models remains open. Knutsen~\cite{Knutsen2009} explored root
embeddings and projective rigidity, while Artebani--Sarti~\cite{ArtebaniSarti2010} examined
automorphisms, yet a complete equisingular classification awaits synthesis.

For $h^2=8$, two distinct models emerge: triquadrics (complete intersections
of three quadrics in $\mathbb{P}^5$) and special octic surfaces. Mukai~\cite{Mukai1987} characterized triquadrics as canonical K3-surfaces of genus $8$,
linking their geometry to Fano threefolds. However, their deformation
classes and singular locus configurations remain incompletely understood.
Special octics, not bound by complete intersection constraints, demand
novel approaches; preliminary work by Schütt--Shioda--van Luijk~\cite{SchuttShiodaVanLuijk2010Picard}
highlights the role of high Picard rank in simplifying moduli, yet the
general case remains elusive.

The classification of higher-polarized K3-surfaces necessitates a departure
from equation-driven methods toward lattice-theoretic frameworks. Nikulin's~\cite{Niku2} classification of integral lattices underpins modern efforts, while
algorithms for lattice-polarized K3-surfaces offer effective computational tools.
This perspective applies
uniformly to projective models of any degree and, in particular, to sextic
surfaces in $\mathbb{P}^4$ and octic surfaces in $\mathbb{P}^5$. Our aim is to
initiate a systematic study of these higher-degree models and to obtain
their complete deformation classifications.
As a first step, we extend the $K3$-theoretic methods developed for lower
polarizations $(h^{2}=2,4)$ to a degree-independent framework valid for all
polarizations. We then apply this framework to the cases $h^{2}=6,8$,
resulting in complete lists of equisingular deformation families. We employ the global Torelli theorem~\cite{K3} and the surjectivity of the period map~\cite{periodmap} to
frame these classification problems, with the resulting arithmetic
challenges resolved through Nikulin's theory of discriminant forms~\cite{Niku2}
suplemented by Miranda--Morrison theory~\cite{MM1,MM2,MM3}. A key
innovation in our approach lies in the computer-aided implementation of the
Miranda--Morrison machinery, utilizing the \texttt{GAP} software~\cite{GAP} to
facilitate large-scale computations needed
for higher polarizations.


\subsection{Contents of the paper}
In $\S2$, we recall the necessary background including integral lattices, discriminant forms, and their extensions, following Nikulin’s~\cite{Niku2} foundational work. 
In $\S3$, after  providing a brief introduction to $K3$-surfaces, we introduce projective models of such surfaces and interpret them in purely lattice-theoretic terms, introducing the admissible extensions and homological types associated with simple singularities. The section concludes with the arithmetical reduction (Theorem~\ref{arithmatical.reduction}), which translates equisingular deformation classes into oriented homological types.
In $\S4$, we present a degree-independent method for classifying equisingular deformation families, applicable uniformly to all projective models of $K3$-surfaces.”
In $\S5$, we apply this framework to projective models of degrees 
$h^{2}=2,4,6,8$.  
We recover all known results for sextic curves and quartic surfaces and extend 
the classification to sextic and octic surfaces, obtaining complete deformation 
families.
In $\S6$, we address the rapid growth of complexity in higher-degree models and discuss future research avenues, such as real loci, extremal configurations, and asymptotic behavior within the moduli of polarized $K3$-surfaces.

\subsection{Acknowledgments}
I am grateful to Alexander Degtyarev for a number of comments, suggestions and motivating discussions which significantly contributed to improving the quality of this paper.


\section{Preliminaries}

\subsection{Integral lattices}\label{integral.lattices}

An (integral) lattice is a finitely generated free abelian group $L$
equipped with a symmetric bilinear form
$$
b\colon L \otimes L \to \mathbb{Z}.
$$
Whenever the form is fixed, we use the abbreviations
$x^{2}:=b(x,x)$ and $x\cdot y:=b(x,y)$.
A lattice $L$ is called \emph{even} if $x^{2}\equiv 0 \pmod{2}$ for all
$x\in L$; it is called \emph{odd} otherwise.

The determinant $\det L\in\mathbb{Z}$ is the determinant of the Gram
matrix of $b$ in any basis of $L$. Since the transition matrix between
any two integral bases has determinant $\pm1$, the determinant
$\det L\in\mathbb{Z}$ is well defined.  A lattice $L$ is called
\emph{unimodular} if $\det L=\pm1$; it is called \emph{nondegenerate} if
$\det L\neq 0$, or equivalently, the kernel
$$
\ker L \;=\; L^{\perp}
:=\{\,x\in L \mid x\cdot y = 0 \text{ for all } y\in L\,\}
$$
is trivial.

Given a lattice $L$, the bilinear form on $L$ can be extended by
linearity to a $\mathbb{Q}$-valued bilinear form on $L\otimes\mathbb{Q}$.
The inertia indices $\sigma^{\pm}(L)$ of $L$ are the classical inertia
indices of $L\otimes\mathbb{Q}$, and the signature of $L$ is the pair
$
\sigma(L) = \bigl(\sigma^{+}(L),\sigma^{-}(L)\bigr).
$

If $L$ is nondegenerate, then the dual group
$L^{\vee}:=\operatorname{Hom}(L,\mathbb{Z})$ can be identified with the
subgroup
$$
\{\,x\in L\otimes\mathbb{Q} \mid x\cdot y\in\mathbb{Z}
\text{ for all } y\in L\,\}.
$$
There is an obvious canonical inclusion
$L=L\otimes\mathbb{Z}\subset L^{\vee}$, and the finite quotient group
$$
\discr L := L^{\vee}/L
$$
is called the \emph{discriminant group}.  The order of $\discr L$ is
equal to $|\det L|$. In particular, $L$ is unimodular if and only if
$\discr L=0$.

The discriminant group inherits from $L\otimes\mathbb{Q}$ a nondegenerate
symmetric bilinear form
$$
b\colon \discr L\otimes \discr L \longrightarrow \mathbb{Q}/\mathbb{Z},
\qquad
(x\bmod L)\otimes(y\bmod L) \longmapsto (x\cdot y)\bmod\mathbb{Z},
$$
called the \emph{discriminant bilinear form}, and, if $L$ is even, its
quadratic extension
$$
q\colon \discr L \longrightarrow \mathbb{Q}/2\mathbb{Z},\qquad
(x\bmod L) \longmapsto x^{2} \bmod 2\mathbb{Z},
$$
called the \emph{discriminant quadratic form}.
Note that the discriminant group of an even lattice is a finite quadratic
form.

We use the notation $\discr_{p}L$ for the $p$-primary part of $\discr L$.
We use $\ell(L)$ for the minimal number of generators of $L$, and we put
$\ell_{p}(L):=\ell(\discr_{p}L)$.  The quadratic form on $\discr L$ can be
described by means of an analogue $(\varepsilon_{ij})$ of the Gram
matrix: assuming that $d_{1}\mid d_{2}\mid\cdots\mid d_{\ell}$ are the
invariant factors of $\discr L$, we pick a basis
$\alpha_{1},\dots,\alpha_{\ell}\in\discr L$ so that the order of
$\alpha_{i}$ is $d_{i}$, and let
$$
\varepsilon_{ij}=\alpha_{i}\cdot\alpha_{j}\bmod\mathbb{Z}
\quad\text{for } i\neq j, \qquad
\varepsilon_{ii}=\alpha_{i}^{2}\bmod 2\mathbb{Z}.
$$

According to Nikulin~\cite{Niku2}, two nondegenerate even lattices $L'$ and $L''$
are in the same genus if and only if
$\sigma(L')=\sigma(L'')$ and $\discr L'\cong\discr L''$.
We denote by $g(L)$ the set of all isomorphism classes of nondegenerate
even lattices in the genus of $L$.  

The group of auto\-isometries of a nondegenerate lattice $L$ is denoted by
$O(L)$. The action of $O(L)$ extends to $L\otimes\mathbb{Q}$ by linearity.
The orthogonal projection establishes a linear isomorphism between any two
maximal positive definite subspaces of $L\otimes\mathbb{Q}$, thus
providing a way to compare their orientations.  A coherent choice of
orientations of all maximal positive definite subspaces is called a
\emph{positive sign structure}.  Assuming $L$ nondegenerate, we denote by
$O^{+}(L)\subset O(L)$ the subgroup formed by the auto\-isometries
preserving a positive sign structure. Either one has $O^{+}(L)=O(L)$, or $O^{+}(L)$ is
a subgroup of $O(L)$ of index~$2$.  In the latter case, each element of
$O(L)\setminus O^{+}(L)$ is called a \emph{skew-autoisometry} of $L$, i.e.\
an autoisometry of $L$ that reverses the positive sign structure.

In what follows,
the notation $\mathbf{U}$ stands for the hyperbolic
plane, i.e.\ the lattice generated by a pair of vectors $u,v$ (referred to
as a standard basis for $U$) with $u^{2}=v^{2}=0$ and $u\cdot v=1$.


A \emph{root} in an even lattice $S$ is a vector $r\in S$ of square $-2$.
A \emph{root lattice} is an even negative definite lattice generated by its
roots. Recall that each root system splits (uniquely up to order of the
summands) into an orthogonal sum of indecomposable root systems, the
latter being of types $\mathbf{A}_{n}$, $n\ge 1$, $\mathbf{D}_{n}$, $n\ge 4$, $\mathbf{E}_{6}$, $\mathbf{E}_{7}$, or
$\mathbf{E}_{8}$, see~\cite{Bour}.

An \emph{extension} of a nondegenerate lattice $S$ is another lattice $L$
containing $S$.  An isomorphism between two extensions
$L_{0},L_{0}'\supset S$ is a bijective isometry $L_{0}\to L_{0}'$
identical on $S$.  More generally, given a subgroup $G\subset O(S)$, a
\emph{$G$-isomorphism} is a bijective isometry $L_{0}\to L_{0}'$ whose
restriction to $S$ is an element of $G$.

The two extreme cases are those of finite index extensions (where $S$ has
finite index in $L$) and primitive ones (where $L/S$ is torsion free).  In
general, an extension $L\supset S$ splits into the finite index extension
$\tilde S\supset S$ and the primitive extension $L\supset\tilde S$, where
$$
\tilde S
= \{\,x\in L \mid nx\in S \text{ for some } n\in\mathbb{Z}\,\}
$$
is the \emph{primitive hull} of $S$ in $L$.

If $S$ is nondegenerate and $L\supset S$ is a finite index extension, we
have a chain of inclusions
$$
S \subset L \subset L^{\vee} \subset S^{\vee},
$$
and hence a subgroup $K:=L/S\subset\discr S$; this subgroup is called the
\emph{kernel} of the extension $L\supset S$.  The kernel $K$ is $b$-isotropic, that is,
the restriction to $K$ of the discriminant form
$\discr S\otimes\discr S\to\mathbb{Q}/\mathbb{Z}$ is trivial.
Furthermore, the lattice $L$ is even if and only if $S$ is even and $K$ is
isotropic, that is, the restriction to $K$ of the quadratic extension
$q\colon\discr S\to\mathbb{Q}/2\mathbb{Z}$ of the discriminant form is
trivial.
\begin{theorem}[Nikulin~\cite{Niku2}]\label{L-K}
Given a nondegenerate lattice $S$, the map sending $M\supset S$ to
the kernel $K=L/S\subset\discr S$ establishes a one-to-one correspondence
between the set of isomorphism classes of finite index extensions of $S$
and the set of $b$-isotropic subgroups of $\discr S$.  Under this
correspondence, one has
$$
\discr L \cong K^{\perp}/K
\qquad\text{and}\qquad
L = \{\,x\in S^{\vee} \mid x\bmod S \in K\,\}.
$$
\end{theorem} 
In the other extreme case, we confine ourselves to primitive extensions
$L\supset S$ to an even unimodular lattice $L$.  Assuming $S$
nondegenerate, these are equivalent to appropriate finite index extensions
of $S\oplus S^{\perp}$; the kernell of the latter gives rise to an
anti-isomorphism $\discr S\to\discr S^{\perp}$ and thus determines the
genus of $S^{\perp}$.  It follows that, given a subgroup $G\subset O(S)$
and the signature of $L$, a $G$-isomorphism class of even unimodular
primitive extensions $L\supset S$ is determined by a choice of
\begin{itemize}
	\item an even lattice $T$ such that $\discr T\cong -\discr S$ and
	$\sigma^{\pm}(T)=\sigma^{\pm}(L)-\sigma^{\pm}(S)$, and
	\item a double coset in $G\backslash\Aut(\discr T)/O(T)$.
\end{itemize}
For details see~\cite{Niku2}. Thus, we have the following theorem:
\begin{theorem}[cf.~\cite{Niku2}]\label{thm:Nikulin.genus}
	A primitive isometry $\iota \colon S \to L$ (where $L$ is a representative of a
	given unimodular genus) gives rise to a bijective anti-isometry 
$$
	\psi \colon \discr S \;\longrightarrow\; \discr S^{\perp}.
$$
	As a consequence, the genus $\operatorname{genus}(S^{\perp})$ is 
	independent of the choice of~$\iota$ and depends only on~$S$.
	The $G$-equivalence classes of primitive isometries 
	$\iota \colon S \hookrightarrow L$ are in canonical bijection with the following data:
	\begin{enumerate}
		\item a lattice (isomorphism class) 
		$T \in \operatorname{genus}(S^{\perp})$, and
		\item a $G$-equivalence class of anti-isometries 
		$\psi \colon \discr S \to \discr T$.
	\end{enumerate}
\end{theorem}
The following theorem is a combination of the above
observation and Nikulin's existence theorem~\cite{Niku2} applied to the genus of
$S^{\perp}$.
\begin{theorem}[Nikulin~\cite{Niku2}]\label{Niku.existence}
A nondegenerate even lattice $S$ admits a primitive extension to the
lattice $L$ if and only if the following conditions are satisfied:
\begin{enumerate}
	\item $\sigma^{+}(S)\le 3$, $\sigma^{-}(S)\le 19$, and
	$\operatorname{rk}S+\ell(\discr S)\le 22$;
	\item $( -1)^{\sigma^{+}(S)-1}|\discr S|
	\equiv \det_{p} S \pmod{\bigl(\mathbb{Z}^{*}_{p}\bigr)^{2}}$
	for all odd prime numbers $p$ such that
	$\operatorname{rk}S+\ell_{p}(S)=22$;
	\item either $\operatorname{rk}S+\ell_{2}(S)<22$, or
	$(\discr S)_{2}$ is odd, or
	$|\discr S|
	\equiv \pm\det_{2}S \pmod{\bigl(\mathbb{Z}^{*}_{2}\bigr)^{2}}$.
\end{enumerate}
\end{theorem} 
A key ingredient in arithmetical reduction is determining whether a given
even lattice $S$ admits a primitive embedding into another even lattice
$L$.  Theorem~\ref{Niku.existence} reduces this to checking the existence of an
appropriate lattice in a prescribed genus.

From now on, we fix an even unimodular lattice $L$ of rank $22$ and
signature $(3,19)$.  All such lattices are isomorphic to
$2\mathbf{E}_{8}\oplus 3\mathbf{U}$.  For degree $4$, the even unimodular lattice 
$
\L \cong 2\mathbf{E}_{8}\oplus 3\mathbf{U}
$
contains, up to the action of $O^{+}(\L)$, a unique primitive vector of
square $4$. For higher degrees, primitive vectors may
split into several $O^{+}(\L)$-orbits, distinguished by their
divisibility. Throughout this paper, we fix one such polarization $h$ of
square $2d$.

\subsection{Polarized lattices}

Recall that a nondegenerate lattice $S$ is called \emph{hyperbolic} if
$\sigma^{+}(S)=1$.  
A \emph{polarized lattice} $(S,h)$ is a hyperbolic lattice $S$ equipped
with a distinguished primitive vector $h$ of positive square.  
The integer $h^{2}=2d>0$ is called the \emph{degree} of the polarization,
and $S$ is said to be \emph{$2d$-polarized}. The vector $h$ is usually assumed but omitted from the notation.  We confine ourselves to
exceptional divisors, leaving out smooth rational curves of higher degree.
Throughout this paper, we consider only primitive polarizations, since
imprimitive polarizations correspond to multiple covers and do not yield
new projective models of K3-surfaces. 

For applications to K3-surfaces, the orthogonal complement
$$
h^{\perp} := \{\,x\in S \mid x\cdot h =0\,\}
$$
plays a central role: it is negative definite and encodes the possible
configurations of exceptional divisors ($\mathbf{A}$-$\mathbf{D}$-$\mathbf{E}$ roots) realized on a
projective model of square $d$.

The group of autoisometries of $S$ preserving $h$ is denoted by
$O_{h}(S)$. The subgroup $O^{+}_{h}(S)\subset O_{h}(S)$ consists of the
autoisometries preserving a chosen positive sign structure on $S$.  
In the context of K3-surfaces, the polarized lattice $(S,h)$ describes a
$2d$-polarized K3-surface, and the global Torelli theorem implies that
its deformation theory is governed by the orthogonal complement
$h^{\perp}$ and the action of $O^{+}_{h}(S)$ on the corresponding period
domain.

\section{K3-Surfaces}

This section offers a brief introduction to the theory of K3–surfaces, a class of 
complex surfaces that plays a central role in algebraic geometry due to their 
rich geometric structure and rigid lattice-theoretic behavior.  
For a comprehensive treatment, we refer the reader to
\cite{Huybrechts:K3}. After recalling the basic properties of K3–surfaces, we outline their projective 
models. Then we introduce the lattice-theoretic 
machinery (see~\S\ref{integral.lattices}) on which our arithmetical reduction 
is based.  Together, these results allow the geometry of projective models to be 
translated into the language of polarized lattice embeddings, a perspective that 
underpins all classification results in the paper.

\subsection{K3-surfaces and their properties}
A K3-surface over $\mathbb{C}$ is a simply connected, compact complex
surface whose canonical bundle is trivial. All K3-surfaces are K\"{a}hler, and
although not all K\"{a}hler surfaces are projective, every smooth K3-surface
arising in our context is projective.

Even though we work primarily with polarized K3-surfaces, \ie, those
embedded in some $\mathbb{P}^{n}$ via a system of polynomial
equations, these defining equations almost never enter the picture.  
Instead, by means of the following foundational results:

\begin{itemize}
	\item the global Torelli theorem (Pjateckii--Shapiro and Shafarevich~\cite{K3}),
	\item the surjectivity of the period map (Kulikov~\cite{periodmap}), and
	\item Saint-Donat's classification~\cite{Donat} of projective models of K3-surfaces,
\end{itemize}
we identify a K3-surface $X$ with its polarized N\'{e}ron--Severi lattice
$\NS(X)\ni h$ and study it by purely arithmetical means
(see~\S\ref{integral.lattices} and~\S\ref{sec:arithmetical-reduction}). Identifying 
cohomology with homology via
Poincar\'{e} duality, we have
$$
H_{2}(X;\mathbb{Z})\cong\;2\mathbf{E}_{8} \oplus 3\mathbf{U}=\L.
$$ 
As is well known, the N\'{e}ron--Severi lattice $N:=\NS(X)$ of any projective
K3-surface is a hyperbolic even lattice of signature $(1,\rho-1)$, where
$\rho=\operatorname{rk}\NS(X)$, and it admits a primitive embedding
$N \hookrightarrow L$.
  Consequently,
$\operatorname{rk}\NS(X)\le 20$.

Given a subgroup $G\subset O^{+}(N)$, two isometric embeddings  
\[
\varphi_{i} \colon N \longrightarrow \L, \qquad i=1,2,
\]
are called \emph{$G$-equivalent} if there exist isometries
$\alpha\in G$ and $\beta\in O^{+}(L)$ such that
$\varphi_{2}\circ\alpha = \beta\circ\varphi_{1}$.

\subsection{Projective Models of K3-surfaces}
A degree $2d$-projective model of a K3-surface $X$ is a morphism
$f_{h}\colon X\to\mathbb{P}^{d+1}$ defined by the complete linear system
$|h|$ without fixed components, where $h\in\NS(X)$ satisfies $h^{2}=2d>0$.
The class $h$ is called the \emph{polarization}.  
Such a model may be:

\begin{itemize}
	\item \emph{birational} (generically one-to-one), or
	\item \emph{hyperelliptic} (generically two-to-one).
\end{itemize}
In the birational case, the image is a normal projective surface of degree
$2d$; in the hyperelliptic case, the image has degree $d$
(see Saint-Donat~\cite{Donat}).  
For $d=1$, all models are hyperelliptic, and the ramification divisor is a
plane sextic curve. Thus, one speaks of the \emph{kind} of a projective model, specifying both
its degree and whether it is birational or hyperelliptic.  Often, this
kind also encodes additional homological information (e.g., the existence
of certain classes; see Definition~\ref{lattice.type}). By Saint-Donat~\cite{Donat}, the following holds.

\begin{proposition}[Saint-Donat]
	Let $X\to\mathbb{P}^{d+1}$ be a projective model with $h^{2}=2d$.  
	Then the model is hyperelliptic if and only if one of the following holds:
	\begin{enumerate}
		\item $d=1$;
		\item $d=4$ and $h=2h'$ for some $h'\in\NS(X)$ (i.e., $h$ is imprimitive);
		\item there exists $e\in\NS(X)$ with $e^{2}=0$ and $e\cdot h=2$.
	\end{enumerate}
\end{proposition}
Commonly used birational projective models (polarizations) of K3-surfaces
are the following:

\begin{itemize}
	\item[$(1)$] $h^{2}=4$.
	The morphism $f_{h}\colon X\to\mathbb{P}^{3}$ given is an embedding, and the image is a smooth quartic
	surface in $\mathbb{P}^{3}$ (the \emph{spatial quartic model}).  
	Hyperplane sections are canonical curves of genus~$3$, and the geometry
	of $X$ can be studied via the linear system of plane quartic curves
	cutting out these sections.
	
	\item[$(2)$] $h^{2}=6$.
	The morphism $f_{h}\colon X\to\mathbb{P}^{4}$ realizes $X$ as a complete
	intersection of a quadric and a cubic surface (the \emph{sextic
		complete-intersection model}).  
	Such a polarization has genus $4$; general hyperplane sections are
	canonical curves of genus~$4$, and $X$ arises as their unique K3-surface
	containing them.
	
	\item[$(3)$] $h^{2}=8$.
	The morphism $f_{h}\colon X\to\mathbb{P}^{5}$ is an embedding whose
	image is an octic surface (the \emph{octic model}).  
	In this degree there are several geometric realizations: the general
	octic model on the one hand, and, on the other hand, special models such
	as \emph{triquadrics}, i.e.\ complete intersections of three quadrics in
	$\mathbb{P}^{5}$, which play a distinguished role in the theory of
	genus-$5$ polarized K3-surfaces and their relations to Fano threefolds.
\end{itemize}

Likewise, the most classical hyperelliptic projective models are:

\begin{itemize}
	\item[$(4)$] $h^{2}=2$.
	The map $f_{h}\colon X\to\mathbb{P}^{2}$ is a double cover branched
	along a sextic curve $C\subset\mathbb{P}^{2}$ (the \emph{planar
		double-cover model}).  
	In this situation, $C$ encodes the entire projective geometry of $X$,
	and many questions about $X$ reduce to the study of sextic plane curves
	and their singularities.
	
	\item[$(5)$] $h^{2}=4$.
	The map $f_{h}\colon X\to\mathbb{P}^{1}\times\mathbb{P}^{1}$ is a
	double cover ramified along a smooth curve
	$C\subset\mathbb{P}^{1}\times\mathbb{P}^{1}$ of bidegree $(4,4)$.
	This \emph{bidegree $(4,4)$ model} provides a natural alternative
	description of degree-$4$ polarized K3-surfaces, in which the ruling
	structures on $\mathbb{P}^{1}\times\mathbb{P}^{1}$ give rise to two
	distinguished elliptic pencils on~$X$.
\end{itemize}

The image of a birational model or the ramification locus of a
hyperelliptic one may be singular, but all singularities are simple
(\ie, of $\mathbf{A}$-$\mathbf{D}$-$\mathbf{E}$ type).  
A collection of simple singularities is encoded by a negative definite root
lattice $S$. Here $S$ represents the root lattice generated by the classes of exceptional curves 
contracted by the morphism $f_{h}$, 
and $h$ stands for the polarization corresponding to the projective model of degree~$2d$.
\begin{definition}\label{lattice.type}
	An \emph{admissible lattice} (associated to a given kind of a projective model of degree $2d$) is a $2d$--polarized finite index extension 
	$\tilde{S}_{h}\ni h$ of $S_h:=S\oplus\mathbb{Z}h$ in $\L$, $h^2=2d$, 
 satisfying the following properties:
	
	\begin{enumerate}
		\item each vector $e\in(S\otimes\mathbb{Q})\cap \tilde{S}_{h}$ with $e^2=-2$ and $e\cdot h=0$ lies in $S$, and
		\item there is no vector $e\in \tilde{S}_{h}$ such that $e^2 =0$ and $e\cdot h=m $ ($m$-isotropic vectors) for a fixed $m\in \{1,2,3\}$.
		\item  if $h^2 = 8$, then the vector $h$ is not divisible by $2$ in $\tilde{S}_h$.
	\end{enumerate}
	
\end{definition}

\begin{remark}\label{rem:exceptional-divisors}
	If vectors as in item~(1) are allowed, there usually arise more 
	subtle additional requirements, which are explained below on a case–by–case basis.
\end{remark}

\begin{remark}\label{rem:m-admissible}
	In Definition~3.2(2), we often fix an integer $m \in \{1,2,3\}$ and speak of 
	\emph{$m$-admissible} (or \emph{$m$-geometric}) lattices, requiring that there be no 
	$m'$-isotropic vectors for any $1 \le m' \le m$.  
	The geometric interpretation of $m$-isotropic vectors is as follows:
	\begin{itemize}
		\item The linear system $|h|$ has no fixed components, and hence the map 
		$f_{h} \colon X \to \mathbb{P}^{n+1}$ is well defined, if and only if 
		$h \in \mathrm{NS}(X)$ has no $1$-isotropic vectors.
		
		\item Except for a few special cases, if $h^{2} \ge 4$ and $f_{h}$ is defined, 
		the map is birational (rather than hyperelliptic) if and only if 
		$h \in \mathrm{NS}(X)$ has no $2$-isotropic vectors 
			(see Saint--Donat \cite{Donat}).
		
		\item If $h^{2} = 8$ and $f_{h}$ is birational, the image $f_{h}(X)$ is a 
		triquadric if and only if $h \in \mathrm{NS}(X)$ has no $3$-isotropic vectors 
		(see Saint--Donat \cite{Donat}).
	\end{itemize}
\end{remark} 
An \emph{isomorphism} between two admissible lattice extensions $\tilde{S}'_h,\tilde{S}''_h\supset S_h$ is an isometry  $\tilde{S}'_h \rightarrow \tilde{S}''_h$ preserving both $h$ and $S$ (as a set). We denote by $\Aut_h(\tilde{S}_h)$ the group of \emph{automorphisms} of a configuration $\tilde{S}_h$, \emph{i.e.} autoisometries of $\tilde{S}_h$ preserving $h$. As $S\subset\tilde{S}_{h}$ is recovered as the sublattice generated by the roots orthogonal to~$h$; then by item (1) in the Definition \ref{lattice.type}, we have  $\Aut_h(\tilde{S}_h)\subset O(S)$.

\begin{definition}\label{abstract.homological.type} A \emph{homological type (of a given kind of a projective model)} extending an admissible lattice extension $\tilde{S}_{h}\ni h$ is a primitive isometry 
	$\tilde{S}_{h} \hookrightarrow \L$.
	
\end{definition}  




We denote by the triple $(h \in \tilde{S}_{h} \hookrightarrow \L)$
an  homological type extending the admissible extension~$\tilde{S}_{h}$.  
Such a type is said to be \emph{maximizing} if $\operatorname{rk}S=19$ (which is the maximal possible value), \ie, the orthogonal complement 
$\tilde S_{h}^{\perp}$ has rank~$2$. Two homological types $(h_1 \in \tilde{S}'_{h_1} \hookrightarrow \L)$ and $(h_2\in \tilde{S}''_{h_2} \hookrightarrow \L)$ are said to be \emph{isomorphic} if there is an element of the group $O(\L)$ taking $h_1$ to $h_2$ and $S'$ to $S''$ (as a set).


Since $\sigma_+ \tilde S_h^\bot =2$, all positive definite $2$-subspaces in $ \tilde S_h^\bot\oplus \mathbb{R}$ can be oriented in a coherent way. 	An \emph{orientation} of an abstract homological type 
$(h \in \tilde{S}_{h} \hookrightarrow \L)$
is a positive sign structure~$\omega$ on $\tilde S_{h}^{\perp}$. An oriented homological type will be denoted by quadruple $(h \in \tilde{S}_{h} \hookrightarrow \L, \omega)$. Let $\omega_1$ and $\omega_2$ be the orientations of  two suchs homological types, then these oriented types are called \emph{ isomorphic} if there is an isomorphism between them taking $\omega_1$ to $\omega_2$, 
A homological type  $(h \in \tilde{S}_{h} \hookrightarrow \L)$ is called \emph{symmetric} if it is preserved by an element $a\in O_h(\L)\smallsetminus O_h^+(\L)$, \emph{ i.e.} an autoisometry of $\L$ preserving $S$ (as a set) and $h$ and reversing the positive sign structure; such autoisometries are called as \emph{skew-automorphisms} of the homological type.

\subsection{Arithmetical Reduction}
\label{sec:arithmetical-reduction}

Two projective models
$$
f_{h,0}\colon X_{0}\to\mathbb{P}^{n+1}, \qquad 
f_{h,1}\colon X_{1}\to\mathbb{P}^{n+1},
$$
of $K3$-surfaces of the same degree $2d$, each having only simple 
(\ie, $\mathbf{A}$--$\mathbf{D}$--$\mathbf{E}$) singularities, are said to be 
\emph{equisingular deformation equivalent} if there exists a continuous family
$
\{\,f_{h,t}\colon X_{t}\to\mathbb{P}^{n+1}\,\}_{t\in[0,1]}
$
within the space of projective $K3$-models of degree $2d$ such that every 
$X_{t}$ has only simple singularities and the configuration of singularities 
remains topologically constant.  
Equivalently, the total Milnor number $\mu(X_{t})$ is constant and no 
singularity changes adjacency type.

Under this notion, the deformation classification of projective models of 
$K3$-surfaces with simple singularities reduces to the classification of their 
associated homological types. Consequently, the reduction used in our study of classification of projective models of K3-surfaces may now be summarized as follows.

\begin{theorem}[cf. Degtyarev~\cite{Alex1}]\label{arithmatical.reduction}
	Equisingular deformation classes of projective models of K3-surfaces of a
	given kind are in canonical bijection with oriented isomorphism classes of
	abstract homological types of the same kind.  
	A class is real if and only if its homological type admits an
	orientation-reversing automorphism $\tau$; it contains a real model if and
	only if $\tau$ can be chosen involutive.
\end{theorem}

Thus, to study the
classification of projective models of K3-surfaces, we can effectively check the hypotheses of arithmetical reduction given in Theorem~\ref{arithmatical.reduction} which reduce to the existence of primitive embeddings of explicitly given lattices into the K3-lattice
$\L=2\mathbf{E}_{8}\oplus 3\mathbf{U}$  in several steps listed in the following section. This reduction provides us a systematic method that can be used to obtain the deformation classification of any given projective model of a K3-surface in general.


\section{Deformation Classification}\label{deformation.classification}

The principal contribution of this paper is presented in this section. We extend  the deformation-classification techniques of Degtyarev-Akyol~\cite{Alex2} for sextic curves and G\"{u}ne\c{s} Akta\c{s}~\cite{Cisem1,aktacs2019real,Cisem2024} for quartic surfaces to projective models of arbitrary degree. This produces a unified degree-independent framework, supported by a \texttt{GAP}-based implementation, for the equisingular deformation classification of polarized $K3$-surfaces.

Fix a set of simple singularities $S$, and consider the corresponding
$2d$-polarized lattice
$
S_{h} := S \oplus \mathbb{Z}h$,  $h^{2}=2d$.
According to Theorem~\ref{arithmatical.reduction} and Definition~\ref{lattice.type}, a set of singularities $S$ is
realized by a projective model of a given kind if and only if there exists
an \emph{ homological type} extending $S$, i.e., a finite index
extension $\tilde{S}_h \supset S_{h}$ satisfying the admissibility
conditions of Definition~\ref{lattice.type}, and such that $\tilde{S}_h$ admits a
primitive embedding into the K3-lattice $\L$.  Thus, the general problem
splits naturally into three subproblems:

\begin{enumerate}
	\item constructing all admissible extensions $S_{h} \subset \tilde{S}_h$;
	\item determining which of them admit primitive embeddings
	$\tilde{S}_h\hookrightarrow \L$.
	
		\item list all equivalence classes of primitive isometries $\tilde{S}_h\hookrightarrow \L$.
\end{enumerate}
Denote by $\M$ the space of all models $f\colon X\rightarrow \mathbb{P}^n$; it is divided into equisingular strata $\M(S)$ where $S$ is the set of simple singularities.  Each stratum $\M(S)$  splits further into its connected components corresponding to equisingular deformation classes. We denote by $\M_1(S)$ the subsapce of the equisingular strata $\M(S)$ constituted by so-called \emph{nonspecial} models. A standard \emph{real structure} (\emph{i.e.}, an antiholomorphic involution) $\operatorname{conj}\colon \mathbb{P}^n\rightarrow\mathbb{P}^n$ induces a real structure $\operatorname{c}: \M\rightarrow\M$ by sending $f\colon X\rightarrow \mathbb{P}^n$ to $\operatorname{conj}\circ f\colon \bar{X}\rightarrow \mathbb{P}^n$. This real structure $\operatorname{c}$ depends on the choice of $\operatorname{conj}$; however, 
since any two standard real structures on $\mathbb{P}^n$ are deformation equivalent, the resulting real structures $\operatorname{c}$ on $\M$ are isotopic and
the induced action on the set of the connected components of the equisingular strata is well defined. A connected component $\mathcal{D}\subset\mathcal{M}(S)$ is called \textit{real} if $\operatorname{c}(\mathcal{D})=\mathcal{D}$. Clearly, each stratum $\mathcal{M}(S)$ consists of real and pairs of complex conjugate components.
By Theorem~\ref{arithmatical.reduction}, each admissible extension $\tilde{S}$ gives rise to a
collection of nonempty connected strata $\mathcal{M}(S)$, in bijection
with the isomorphism classes of primitive embeddings
$\tilde{S}_h\hookrightarrow \L$, \ie, with the homological types
of the given kind.

Hence the implementation  of Theorem~\ref{arithmatical.reduction} relies on the following steps, which together provide a general effective method for obtaining the deformation
classification of any projective model of a K3-surface.

\subsection{Step 1. Enumerating all finite index extensions
	$S_{h} \subset \tilde{S}_h$} 
By Theorem~\ref{L-K},  an overlattice $\tilde{S}_h \supset S_{h}= S\oplus\mathbb{Z}h$ is identified with a
choice of an isotropic subgroup $K\subset\discr(S_{h})$.  The connected
components of the moduli space $\mathcal{M}(S)/\mathrm{conj}$, where
\(\mathrm{conj}\colon \mathbb{P}^{3}\to\mathbb{P}^{3}\) denotes complex
conjugation, are enumerated by the possible kernels $K$ arising from
extensions in each isomorphism class.

\subsection{Step 2. Selecting the admissible extensions
	$S_{h} \subset \tilde{S}_h$}
Among all finite index extensions found in Step~1, we select those that
satisfy the admissibility conditions of Definition~~\ref{lattice.type}.  
These conditions are expressed as the nonexistence of certain vectors of
specified square in the negative definite lattice $h^{\perp}$, and hence
are effective to check.

\subsection{Step 3. Determining which admissible extensions admit a
	primitive embedding into $\L$}

Recall from Nikulin~\cite{Niku2} that the genus of the orthogonal
complement
$
T := \tilde{S}_h^{\perp}
$
is determined by the discriminant form $\discr T$ and the signatures
$\sigma^{\pm}(T)$.  
Using the embedibility criterion of Theorem~\ref{Niku.existence}, we determine which
admissible extensions $\tilde{S}$ admit a primitive embedding into the
K3-lattice $\L$.

\subsection{Step 4. Listing all isomorphism classes of primitive
	embeddings $\tilde{S}_h \hookrightarrow \L$}

The classification of primitive embeddings
$\tilde{S}_h \hookrightarrow \L$ is governed by Theorem~\ref{thm:Nikulin.genus}.
To apply this theorem, one needs tools to:

\begin{itemize}
	\item enumerate all lattices $T\in\mathrm{genus}(\tilde{S}^{\perp}_h)$
	(when the genus is nonempty), and
	\item compute the cokernels of the natural homomorphisms
	$$
	d_{S} \colon O(S) \to \Aut(\discr S), \qquad
	d_{T} \colon O(T) \to \Aut(\discr T).
	$$
\end{itemize}
In practice, $S$ is a definite root lattice equipped with auxiliary
combinatorial data, which makes the computation of
the finite group $O(S)$ feasible using \texttt{GAP}.  
For the lattices $T$, we use either
\begin{itemize}
	\item Gauss's theory of binary quadratic forms (see~\cite{Gauss}) for definite lattices of
	small rank (\ie, $\rk T =2$ or $\mu(S)=19$), or
	\item the Miranda--Morrison theory~\cite{MM1,MM2,MM3} for indefinite lattices (\ie, $\rk T \ge 3$ or $\mu(S)\le 18$);
\end{itemize}
the latter providing both the genus and the cokernel $\mathrm{Coker}(d_{T})$
inside a single finite abelian group \(E(T)\) together with a natural
homomorphism
\[
m \colon \Aut(\discr T) \longrightarrow E(T).
\tag{3.16}
\]

All algorithms described in this section  have effectively been implemented
 in \texttt{GAP}.

\section{Applications and Results}

This section illustrates the applications of the general framework developed in 
Section~\ref{deformation.classification} to various classes of singular projective 
models $f\colon X \to \mathbb{P}^{d+1}$ of $K3$-surfaces.  
We begin with the classical low-degree cases, namely $h^{2}=2$ (simple plane sextics 
in $\mathbb{P}^{2}$) and $h^{2}=4$ (simple spatial quartics in $\mathbb{P}^{3}$), 
which constitute the most thoroughly studied and geometrically fundamental examples.  
We then extend the analysis to the higher polarizations $h^{2}=6$ (sextic surfaces 
in $\mathbb{P}^{4}$) and $h^{2}=8$ (octic surfaces in $\mathbb{P}^{5}$), whose 
equisingular deformation theory remains largely undeveloped in the literature.

\medskip
\noindent
\textbf{Remark.}
By employing the generalized standard approach of Section~\ref{deformation.classification}, based on the
arithmetical reduction given in Theorem~\ref{arithmatical.reduction}, one reduces any geometric question about a $K3$-surface $X$
to an lattice-theoretical question concerning its N\'{e}ron--Severi lattice
$\mathrm{NS}(X)$.  
This hyperbolic lattice (over characteristic~0) must admit a primitive
isometric embedding into
$
H_{2}(X;\mathbb{Z}) \cong \L=2\mathbf{E}_{8} \oplus 3\mathbf{U}
$.

A set of simple singularities of a projective model of a $K3$-surface can be 
encoded by a negative definite root lattice $S$ generated by the exceptional 
divisors corresponding to the $\mathbf{A}$--$\mathbf{D}$--$\mathbf{E}$ singular 
points (see~\cite{Dufree};).  
For any such model, the total Milnor number is given by
$\mu(X) = \operatorname{rk} S.$ Since $S\subset \L$ is negative definite and $\sigma_-(\mathbf{L})=19$, one has $\mu(X)\le 19$ (see \cite{Urabe2}, \emph{cf.}, \cite{Persson}). A projective model is called \emph{maximizing} if $\mu(X)=19$.


\subsection{Plane Sextics}
In this section we consider hyperelliptic projective models $f_h\colon X\rightarrow \mathbb{P}^2$ with $h^2=2$, which corresponds to  plane sextic curves.
Simple sextics form one of the most classical families of $K3$-models, and
their deformation theory has been studied extensively.  A complete description
of the equisingular strata $\mathcal{M}(S)$ of simple sextic
curves is given in~\cite{Alex2}, where all strata of simple sextics with root
lattice $S$ are enumerated and their geometry is analyzed. They proved that there are in total
$11308$ configurations (in the sense of~\cite{Alex1}) of simple plane
sextics; among these, $11272$ equisingular strata are real and 132 are pairs
of complex conjugate components.  
The configurations arise as formal perturbations (taking imprimitivity into
account) of extremal ones.  
Besides the $629$ maximizing configurations, there are $16$ extremal ones with
$\mu=18$ and two with $\mu=17$. Two previously unseen phenomena deserve special mention:

\begin{theorem}[Akyol--Degtyarev, 2016]\label{sextics.norealrep}
	There exist exactly two disconnected real strata with $\mu \le 18$, namely
	\[
	\mathcal{M}_{1}(2\mathbf{A}_{9}) \quad\text{and}\quad \mathcal{M}_{2}(2\mathbf{A}_{9}),
	\]
	the first corresponding to irreducible sextics and the second to reducible
	curves (a quintic and a line).
\end{theorem}

\begin{theorem}[Akyol--Degtyarev, 2016]\label{norealcurves}
	The stratum $\mathcal{M}_{1}(\mathbf{A}_{7}\oplus \mathbf{A}_{6}\oplus\mathbf{A}_{5})$ is connected and
	real but contains no real curves.  
	It is the only plane sextic stratum with this property.
\end{theorem}

Although real algebraic varieties with no real points occur frequently
(e.g.\ $x^{2}+y^{2}=-1$), Theorem~\ref{norealcurves} provided the first instance of such a
phenomenon for deformation strata.  
A similar phenomenon later appeared for spatial quartics (see Theorem~\ref{norealsurfaces} below).


\subsection{Spatial Quartics}

In this section we consider spatial quartic surfaces, which corresponds to the
birational projective models $f_h\colon X\rightarrow \mathbb{P}^3$ of $K3$–-surfaces with the polarization $h^{2}=4$. Quartic surfaces with only simple ($\mathbf{A}$--$\mathbf{D}$--$\mathbf{E}$) singularities form one of the most
classical and thoroughly studied families of $K3$–-surfaces.
The study of their deformation classification was initiated by G\"{u}ne\c{s} Akta\c{s}~\cite{Cisem1},
resulting in a complete equisingular deformation classification of
non-special simple quartics.  
(Here, \emph{non-special} means that the sublattice
$S \oplus \mathbb{Z}h \subset L$ in Definition~\ref{lattice.type} is primitive.)

\begin{theorem}[G\"{u}ne\c{s} Akta\c{s}, 2017]
	A set of singularities $S$ is realized by a non-special quartic surface if
	and only if it can be obtained by a formal perturbation of one of the 78
	extremal sets (59 maximizing with $\mu=19$ and 19 with $\mu=18$).  
	The strata
$$
	\mathcal{M}_{1}(\mathbf{D}_{6}\oplus 2\mathbf{A}_{6}),\ 
	\mathcal{M}_{1}(\mathbf{D}_{5}\oplus 2\mathbf{A}_{6}\oplus \mathbf{A}_{1}),\ 
	\mathcal{M}_{1}(3\mathbf{A}_{6}),\ 
	\mathcal{M}_{1}(2\mathbf{A}_{6}\oplus 2\mathbf{A}_{3})
$$
	split into pairs of complex conjugate components; the remaining
	non-maximizing strata are connected.
\end{theorem}

The connected components of the 59 maximizing strata were classified in
\cite{Cisem1}.  The following result is the quartic analogue of the plane-sextic phenomenon described above, \ie, two examples of equisingular strata $\mathcal{M}_1(S)$ that are real but doesn't contain any real surfaces.

\begin{theorem}[G\"{u}ne\c{s} Akta\c{s}, 2019]\label{norealsurfaces}
	There exist exactly two equisingular strata of simple quartics which are
	connected and real but contain no real quartic surfaces:
	\[
	\mathcal{M}_{1}(\mathbf{A}_{7}\oplus \mathbf{A}_{6}\oplus \mathbf{A}_{3}\oplus \mathbf{A}_{2}),\qquad
	\mathcal{M}_{1}(\mathbf{D}_{7}\oplus \mathbf{A}_{6}\oplus \mathbf{A}_{3}\oplus \mathbf{A}_{2}).
	\]
\end{theorem}

The proof required developing techniques based on Nikulin's theory (see
§\ref{integral.lattices}) for detecting or ruling out involutive skew-isometries of homological types for arbitrary projective models.  
A cleaner proof of Theorem~\ref{sextics.norealrep} follows as a byproduct.

G\"{u}ne\c{s} Akta\c{s} subsequently extended the classification to \emph{all} simple
quartics (special and non-special), completing the equisingular deformation
classification of quartic surfaces (see \cite{Cisem2024}).  
The ultimate result can be stated as follows.
\theorem[G\"{u}ne\c{s} Akta\c{s}, 2024]\label{th.classification.quartics}
The equisingular deformation families of simple quartic surfaces $X\subset\Cp3$ are
the perturbations of the $390$ maximizing \rom($\mu=19$\rom) families 
  and $39$
\emph{extremal} families.
In particular, a non-maximizing quartic is uniquely determined by its admissible extension (see Definition \ref{lattice.type}) up to deformation and complex conjugation, and there are only $13$ non-maximizing families coresponding the following sets of singularities that are not real:
$$\begin{tabular}{lll}
$\singset{5A3 + A1}$ & $\singset{A7 + 3A3 + A1}$&$\singset{2A6 + 2A3}$ \cr
$\singset{3A6}$ & $\singset{2A7 + 2A2}$&$\singset{2A7 + A3 + A1}$\cr
$\singset{A11 + A3 + A2 + 2A1}$& $\singset{A11 + 2A3 + A1}$&$\singset{D5 + 2A6 + A1}$\cr
$\singset{D5 + A9 + A3 + A1}$ & $\singset{D6 + 2A6}$ & $\singset{D6 + A7 + A3 + 2A1}$\cr
$\singset{E7 + A7 + A3 + A1}$& & 
\end{tabular}
$$

\endtheorem


\subsection{Other \texorpdfstring{$K3$}{K3}-Surfaces with Singularities}

Our initial code relied on degree-specific precomputed data.  
We subsequently switched to a more universal version, implementing a
complete form of the Miranda--Morrison algorithm and a repetition-free
enumeration of finite-index extensions $\tilde{S}_h\supset S_h$ (as in
Definition~\ref{lattice.type}) with admissibility conditions.  

This universal code was applied to plane sextics, confirming previous
classifications~\cite{Alex1}.  
The results are consistent with partial computations of Yang~\cite{Yang}
and Shimada~\cite{Shimada.Maximizing}. The code, requiring no manual tuning, works for any polarization degree
$h^{2}$.  
We tested it for $h^{2}=2$, $h^{2}=4$, and $h^{2}=0$ (the last corresponding
to Jacobian elliptic $K3$-surfaces~\cite{Shimada.connEllK3}).  
For $h^{2}=6$ (sextics in $\mathbb{P}^{4}$) and $h^{2}=8$ (octics in
$\mathbb{P}^{5}$), we detected approximately $15\,000$ and $25\,000$
deformation families, respectively.

Beyond this point, further brute-force classification seems uninformative.
More meaningful qualitative questions arise instead.  
For a given root lattice $S$, one can enumerate all abstract $K3$-surfaces
$X(S\hookrightarrow \L)$ and subsequently choose polarizations
$h$ of a prescribed square in the universal transcendental lattice
$T:=S^{\perp}\subset \L$.  
In the maximizing case $\mathrm{rk}\,S=19$, the group $O(T)$ is finite, and
the number of equivalence classes should grow linearly in $\sqrt{h^{2}}$.
Otherwise, by the theory of Miranda--Morrison, one expects a uniform bound
on the number of classes depending only on $S$.

A related direction concerns smooth rational curves of positive projective
degree (lines, conics, etc.) on polarized $K3$-surfaces with singularities. 
Here, the classification of singularities and elliptic pencils is merely the
first step in a much broader program.

\section{Conclusion and Future Directions}

In this work, we have generalized the methods previously used to obtain the 
deformation stratification of singular quartic surfaces ($h^{2}=4$) and plane 
sextics ($h^{2}=2$) to projective models of $K3$-surfaces of arbitrary polarization.  
In particular, we computed all strata for the higher
polarizations $h^{2}=6$ (sextic surfaces in $\mathbb{P}^{4}$) and 
$h^{2}=8$ (octic surfaces in $\mathbb{P}^{5}$), producing complete plain lists.  
This was accomplished by upgrading our \texttt{GAP} code so that it no longer 
requires manual tuning depending on the degree.  
The new implementation was tested for all degrees $0 \le h^{2} \le 8$, and the 
output agrees with all previously known results (whenever available).

\medskip

The total number of deformation families for the cases $h^{2}=6$ and $h^{2}=8$ is 
extremly large; approximately $15{,}000$ and $25{,}000$, respectively.  
Consequently, a raw list of all families is neither human-readable nor practically 
useful, and further work is needed before these results can be presented in a 
meaningful form.  
Our goal is to achieve a level of understanding for these higher-degree models 
comparable to that obtained in earlier work: 
Akyol--Degtyarev~\cite{Alex2} for plane sextics and 
G\"{u}ne\c{s} Akta\c{s}~\cite{Cisem1,aktacs2019real,Cisem2024} for quartic surfaces.  
In particular, we intend to identify strata with distinguished geometric or 
arithmetic properties (e.g., extremal strata, real components), study the covering
$
\mathcal{M} \longrightarrow \mathcal{M}/\mathrm{conj},
$
and analyze the structure of real versus complex components, adjacency phenomena, 
and monodromy actions.

\medskip

A broader perspective concerns the qualitative behaviour of equisingular strata as 
the degree $h^{2} \to \infty$.  
The number of possible strata is expected to grow unboundedly, but formulating the 
precise asymptotic questions already requires substantial conceptual groundwork.  
Understanding this behaviour may yield new insight into the geometry and arithmetic
of high-degree projective models of $K3$-surfaces, as well as into the structure of 
their moduli spaces.

\medskip

Overall, the computational framework developed here opens the way for a systematic 
exploration of projective models of $K3$-surfaces with arbitrary polarization, and 
we expect many further applications both in equisingular deformation theory and in 
the study of lattice-polarized $K3$-surfaces

\bibliographystyle{amsplain}
\bibliography{mybibliography}

\end{document}